\documentclass[12pt]{article}
\usepackage{amsmath}
\usepackage{amssymb}
\usepackage{amsbsy}
\usepackage{fancyhdr}
\usepackage{xcolor}
\usepackage{amsfonts}
\newtheorem{lem}{Lemma}[section]%
\newtheorem{theorem}[lem]{Theorem}%
\newtheorem{cor}[lem]{Corollary}%
\newtheorem{prop}[lem]{Proposition}%
\newtheorem{claim}{Claim}
\newtheorem{case}{Case}

\newtheorem{construction}{Construction}

\def\a{\alpha} \def\b{\beta} \def\g{\gamma} \def\d{\delta} 
 \def\s{\sigma}   
 \def\ld{\lambda} 

\def\G{\Gamma}

 \def\lg{\langle} \def\rg{\rangle}
\def\nd{\mathrel{\bigm|\kern-.7em/}}

\def\f{\noindent}

\def\z{\mathbb{Z}}
\def\PSL{\hbox{\rm PSL}}

\def\GL{\hbox{\rm GL}}
\def\PGL{\hbox{\rm PGL}}

\def\Aut{\hbox{\rm Aut\,}}

\def\Cay{\hbox{\rm Cay }}

\def\mod{\hbox{\rm mod }}

\def\F{\hbox{\rm F}}

\def\BiCay{{\rm BiCay}}
\def\B{{\rm BiCay}(H,R,L,S)}
\def\Bb{{\rm BiCay}(H,\emptyset,\emptyset,S)}

\def\H{{\mathcal H}}
\def\GG{{\mathcal G}}
\def\demo{{\bf Proof}\hskip10pt}

\def\mz{{\mathbb Z}}
\def\qed{\hskip10pt $\Box$\vspace{3mm}}

\begin{document}
\title{Cubic edge-transitive bi-$p$-metacirculants
\author
{Yan-Li Qin\ \ \ \ \ \ {Jin-Xin Zhou}\footnote{\small{This work was partially supported by the National
Natural Science Foundation of China (11271012) and the Fundamental
Research Funds for the Central Universities (2015JBM110).}}\\
{\small Mathematics, Beijing Jiaotong University, Beijing 100044, P.R. China}\\
{\small{Email}:\ 15118413@bjtu.edu.cn, jxzhou@bjtu.edu.cn}
}}

\date{}
\maketitle

\begin{abstract}
A graph is said to be a {\em bi-Cayley graph} over a group $H$ if it admits $H$ as a group of automorphisms acting semiregularly on its vertices with two orbits. For a prime $p$, we call a bi-Cayley graph over a metacyclic $p$-group a {\em bi-$p$-metacirculant}. In this paper, the automorphism group of a connected cubic edge-transitive bi-$p$-metacirculant is characterized for an odd prime $p$, and the result reveals that a connected cubic edge-transitive bi-$p$-metacirculant exists only when $p=3$. Using this, a classification is given of connected cubic edge-transitive bi-Cayley graphs over an inner-abelian metacyclic $3$-group. As a result, we construct the first known infinite family of cubic semisymmetric graphs of order twice a $3$-power.

\bigskip

\noindent{\bf Keywords} bi-$p$-metacirculant, symmetric graph, semmisymmetric graph\\
\noindent{\bf 2000 Mathematics subject classification:} 05C25, 20B25.
\end{abstract}

\section{Introduction}

Throughout this paper, groups are assumed to be finite, and graphs are assumed to be finite, connected, simple and undirected.
For a graph $\G$, we denote by $V(\G)$ the set of all vertices of $\G$, by $E(\G)$ the set of all edges of $\G$, by $A(\G)$ the set of all arcs (ordered paries of adjacent vertices) of $\G$, and by $\Aut(\G)$ the full automorphism group of $\G$. For $u, v\in V(\G)$, denote by $\{u, v\}$ the edge incident to $u$ and $v$ in $\G$.
For the group-theoretic and the graph-theoretic terminology not defined here we refer the reader to \cite{Bondy,Wielandt}.

Let $\G$ be a graph. If $\Aut(\G)$ is transitive on $V(\G)$, $E(\G)$ or $A(\G)$, then $\G$ is said to be {\em vertex-transitive}, {\em edge-transitive} or {\em arc-transitive}, respectively.
An arc-transitive graph is also called a symmetric graph.
A graph $\G$ is said to be {\em semisymmetric} if $\G$ has regular valency and is edge- but not vertex-transitive.

Let $G$ be a permutation group on a set $\Omega$ and $\a\in \Omega$.
Denote by $G_{\a}$ the stabilizer of $\a$ in $G$, that is the subgroup of $G$ fixing the point $\a$.
We say that $G$ is {\em semiregualr} on $\Omega $ if $G_{\a}=1$ for every $\a\in \Omega$ and {\em regular} if $G$ is transitive and semiregular.
A graph is said to be a {\em bi-Cayley graph} over a group $H$ if it admits $H$ as a semiregular automorphism group with two orbits (Bi-Cayley graph is sometimes called {\em semi-Cayley graph}).
Note that every bi-Cayley graph admits the following concrete realization. Given a group $H$, let $R$, $L$ and $S$ be subsets of $H$ such that $R^{-1}=R$, $L^{-1}=L$ and $R\cup L$ does not contain the identity element of $H$. The {\em bi-Cayley graph} over $H$ relative to the triple $(R, L, S)$, denoted by BiCay($H, R, L, S$), is the graph having vertex set the union of the right part $H_{0}=\{h_{0}~|~h\in H\}$ and the left part $H_{1}=\{h_{1}~|~h\in H\}$, and edge set the union of the right edges $\{\{h_{0},~g_{0}\}~|~gh^{-1}\in R\}$, the left edges $\{\{h_{1},~g_{1}\}~|~gh^{-1}\in L\}$ and the spokes $\{\{h_{0},~g_{1}\}~|~gh^{-1}\in S\}$.
Let $\G=\BiCay(H, R, L, S)$. For $g\in H$, define a permutation $R(g)$ on the vertices of $\G$ by the rule
$$h_i^{R(g)}=(hg)_i, \forall i\in\mz_2, h\in H.$$
Then $R(H)=\{R(g)\ |\ g\in H\}$ is a semiregular subgroup of $\Aut(\G)$ which is isomorphic to $H$ and has $H_0$ and $H_1$ as its two orbits.
When $R(H)$ is normal in $\Aut(\G)$, the bi-Cayley graph $\G=\BiCay(H,R,L,S)$ will be called a {\em normal bi-Cayley graph} over $H$ (see \cite{Zhouaut}).

A natural problem in the study of bi-Cayley graphs is: for a given finite group $H$, to classify bi-Cayley graphs with specific symmetry properties over $H$. Some partial answers for this problem have been obtained. For example, in \cite{Boben} Boben et al. studied some properties of cubic $2$-type bi-Cayley graphs over cyclic groups and the configurations arising from these graphs, in \cite{Pisanski} Pisanski classified cubic bi-Cayley graphs over cyclic groups, in \cite{BiCayley.2} Kov\'{a}cs et al. gave a classification of arc-transitive one-matching abelian bi-Cayley graphs,
and more recently, Zhou et al. \cite{Zhoucubic} gave a classification of cubic vertex-transitive abelian bi-Cayley graphs.
In this paper, we shall investigate cubic edge-transitive bi-Cayley graphs over metacyclic $p$-groups where $p$ is an odd prime. Following up \cite{NET}, we call a bi-Cayley graph over a metacyclic $p$-group a {\em bi-$p$-metacirculant}.

Another motivation for us to consider bi-Cayley graphs over metacyclic $p$-groups is the observation that the Gray graph \cite{Bouwer1}, the smallest trivalent semmisymmetric graph, is a bi-Cayley graph over a non-abelian metacyclic group of order $27$. In \cite{NET}, the cubic edge-transitive bi-Cayley graphs over abelian groups have been classified. So, we shall restrict our attention to bi-Cayley graphs over non-abelian metacyclic $p$-groups.

Our first result characterizes the automorphism groups of cubic edge-transitive bi-$p$-metacirculants.

\begin{theorem}\label{5no}
Let $\G$ be a connected cubic edge-transitive bi-Cayley graph over a non-abelian metacyclic $p$-group $H$ with $p$ an odd prime.
Then $p=3$, and either $\G$ is isomorphic to the Gray graph or $\G$ is a normal bi-Cayley graph over $H$.
\end{theorem}

Applying the above theorem, our second result gives a classification of connected cubic edge-transitive bi-Cayley graphs over a inner-abelian metacyclic $p$-group. A non-abelian group is called an {\em inner-abelian group} if  all of its proper subgroups are abelian.

\begin{theorem}\label{result}
Let $\G$ be a connected cubic edge-transitive bi-Cayley graph over an inner-abelian metacyclic $3$-group $H$. Then $\G$ is isomorphic to either $\G_{t}$ or $\Sigma_{t}$ (see Section~\ref{sec-5} for the construction of these two families of graphs).
\end{theorem}

Theorem~\ref{5no} also enables us to give a short proof of the main result in \cite{2p3}.

\begin{cor}{\rm~\cite[Theorem 1.1]{2p3}}\label{cor2p3}
Let $p$ be a prime. Then, with the exception of the Gray graph on $54$ vertices, every cubic edge-transitive graph of order $2p^3$ is vertex-transitive.
\end{cor}


\section{Preliminaries}

In this section, we first introduce the notation used in this paper. For a positive integer $n$, denote by $\z_n$ the cyclic group of order $n$ and
by $\z_n^*$ the multiplicative group of $\z_n$ consisting of numbers coprime to $n$.
For a finite group $G$, the full automorphism group, the center, the derived subgroup and the Frattini subgroup of $G$ will be denoted by $\Aut(G)$, $Z(G)$, $G'$ and $\Phi(G)$, respectively. For $x, y\in G$, denote by $[x, y]$ the commutator $x^{-1}y^{-1}xy$. For a subgroup $H$ of $G$,
denote by $C_G(H)$ the centralizer of $H$ in $G$ and by $N_G(H)$ the
normalizer of $H$ in $G$. For two groups $M$ and $N$, $N\rtimes M$ denotes a semidirect product of $N$ by $M$.

Below, we restate some group-theoretic results, of which the first
is usually called the $N/C$-theorem.

\begin{prop} {\rm~\cite[Chapter 1, Theorem 4.5]{Huppert}}\label{NC}
Let $H$ be a subgroup of a group $G$. Then $C_G(H)$ is normal in $N_G(H)$, and
the quotient group $N_G(H)/C_G(H)$ is isomorphic to a subgroup of $\Aut(H)$.
\end{prop}



Now we give two results regarding metacyclic $p$-groups.

\begin{prop}{\rm~\cite[Lemma~2.4]{pcomplement}}\label{metap}
Let $P$ be a split metacyclic $p$-group:
\begin{center}
$P=\lg x,y\mid x^{p^m}=y^{p^n}=1, yxy^{-1}=x^{1+p^l}\rg$, where $0<l<m$, $m-l\leq n$.
\end{center}
Then the automorphism group $\Aut(P)$ of $P$ is a semidirect product of a normal $p$-subgroup and the cyclic subgroup $\lg\sigma\rg$ of order $p-1$,
where $\sigma(x)=x^r$ and $\sigma(y)=y$, $r$ is a primitive $(p-1)$th root of unity modulo $p^m$.
\end{prop}

\begin{prop}{\rm~\cite[{Proposition~2.3}]{pcomplement}}\label{pcom}
Let $G$ be a finite group with a non-abelian metacyclic Sylow $p$-subgroup $P$. If $P$ is nonsplit, then $G$ has a normal $p$-complement.
\end{prop}

Next, we give some results about graphs. Let $\G$ be a connected graph with an edge-transitive group $G$ of automorphisms and let $N$ be a normal subgroup of $G$. The {\em quotient graph} $\G_N$ of $\G$ relative to $N$ is defined as the graph
with vertices the orbits of $N$ on $V(\G)$ and with two orbits
adjacent if there exists an edge in $\G$ between the vertices lying in
those two orbits. Below we introduce two propositions, of which the
first is a special case of \cite[Theorem~9]{VTgraph}.

\begin{prop}\label{3orbits}
Let $\G$ be a cubic graph and let $G\leq \Aut(\G)$ be arc-transitive on $\G$.
Then $G$ is an $s$-arc-regular subgroup of $\Aut(\G)$ for some integer $s$.
If $N\unlhd G$ has more than two orbits in $V(\G)$, then $N$ is semiregular on $V(\G)$,
$\G_N$ is a cubic symmetric graph with $G/N$ as an $s$-arc-regular subgroup of automorphisms.
\end{prop}

The next proposition is a special case of \cite[Lemma~3.2]{6p2}.

\begin{prop}\label{intransitive}
Let $\G$ be a cubic graph and let $G\leq \Aut(\G)$ be transitive on $E(\G)$ but intransitive on $V(\G)$.
Then $\G$ is a bipartite graph with two partition sets, say $V_0$ and $V_1$.
If $N \trianglelefteq G$ is intransitive on each of $V_0$ and $V_1$, then $N$ is semiregular on $V(\G)$,
$\G_N$ is a cubic graph with $G/N$ as an edge- but not vertex-transitive group of automorphisms.
\end{prop}

The next proposition is basic for bi-Cayley graphs.

\begin{prop}{\rm~\cite[Lemma~3.1]{Zhoucubic}}\label{bicayley}
Let $\G=\BiCay(H,R,L,S)$ be a connected bi-Cayley graph over a group $H$. Then the following hold:
\begin{enumerate}
\item[$(1)$] $H$ is generated by $R\cup L\cup S$.
\item[$(2)$] Up to graph isomorphism, $S$ can be chosen to contain the identity of $H$.
\item[$(3)$] For any automorphism $\alpha$ of $H$, $\BiCay(H, R, L, S)\cong \BiCay(H, R^{\alpha}, L^{\alpha}, S^{\alpha})$.
\item[$(4)$] $\BiCay(H, R, L, S) \cong \BiCay(H, L, R, S^{-1})$.
\end{enumerate}
\end{prop}

Next, we collect several results about the automorphisms of the bi-Cayley graph $\G={\rm BiCay}(H, R, L, S)$.
Recall that for each $g\in H$, $R(g)$ is a permutation on $V(\G)$ defined by the rule
\begin{equation}\label{1}
h_{i}^{R(g)}=(hg)_{i},~~~\forall i\in \mathbb{Z}_{2},~h,~g\in H,
\end{equation}
and $R(H)=\{R(g)\ |\ g\in H\}\leq\Aut(\G)$.
For an automorphism $\a$ of $H$ and $x,y,g\in H$, define two permutations on $V(\G)=H_0\cup H_1$ as following:
\begin{equation}\label{2}
\begin{array}{ll}
\d_{\a,x,y}:& h_0\mapsto (xh^{\a})_1, ~h_1\mapsto (yh^{\a})_0, ~\forall h\in H,\\
\s_{\a,g}:& h_0\mapsto (h^\a)_0, ~h_1\mapsto (gh^{\a})_1, ~\forall h\in H.\\
\end{array}
\end{equation}
Set \begin{equation}\label{3}
\begin{array}{lll}
{\rm I}&=& \{\d_{\a,x,y}\ |\ \a\in\Aut(H)\ s.t.\ R^\a=x^{-1}Lx, ~L^\a=y^{-1}Ry, ~S^\a=y^{-1}S^{-1}x\},\\
{\rm F} &=&\{ \s_{\a,g}\ |\ \a\in\Aut(H)\ s.t.\ R^\a=R, ~L^\a=g^{-1}Lg, ~S^\a=g^{-1}S\}.
\end{array}
\end{equation}

\begin{prop}{\rm~\cite[Theorem~3.4]{Zhouaut}}\label{bicayleyaut}
Let $\Gamma=\BiCay(H, R, L, S)$ be a connected bi-Cayley graph over the group $H$.
Then $N_{\Aut(\Gamma)}(R(H))=R(H)\rtimes F$ if $I=\emptyset$ and
$N_{\Aut(\Gamma)}(R(H))=R(H)\langle F, \delta_{\a,x,y}\rangle$ if
$I\neq\emptyset$ and $\delta_{\a,x,y}\in I$. Furthermore, for any $\delta_{\a,x,y}\in I$, we have the following:
\begin{enumerate}
\item[$(1)$] $\langle R(H), \delta_{\a,x,y}\rangle$ acts transitively on $V(\Gamma)$;
\item[$(2)$] if $\a$ has order $2$ and $x=y=1$, then $\Gamma$ is isomorphic to the Cayley graph $\Cay(\bar{H},~R\cup \alpha S)$, where $\bar{H}=H\rtimes \langle\a\rangle$.
\end{enumerate}
\end{prop}

\begin{prop}{\rm~\cite[Proposition~5.2]{NET}}\label{abelian-biCayley}
Let $n, m$ be two positive integers such that $nm^2\geq 3$. Let $\ld=0$ if $n=1$, and let $\ld\in\z_n^{*}$ be such that $\ld^2-\ld+1\equiv 0\ ({\rm mod}\ n)$ if $n> 1$.
Let
\[
\begin{array}{ll}
\mbox{$\mathcal{H}_{m,n}=\lg x\rg\times\lg y\rg\cong\z_{nm}\times\z_m$,} & \\
\mbox{$\G_{m,n,\ld}={\rm BiCay}(\mathcal{H}_{m,n}, \emptyset, \emptyset, \{1,x,x^{\ld}y\})$.} &
\end{array}
\]
Let $\G=\B$ be a connected cubic normal edge-transitive bi-Cayley graph over an abelian group $H$. Then $\G\cong\G_{m,n,\ld}$ for some integers $m,n,\ld$.
\end{prop}

Finally, we give some results about cubic edge-transitive graphs.

\begin{prop}{\rm~\cite[Theorem~3.2]{2pn}}\label{Sym-P-normal}
Let $\G$ be a connected cubic symmetric graph of order $2p^n$ with $p$ an odd prime and $n$ a positive integer.
If $p\neq 5,7$, then every Sylow $p$-subgroup of $\Aut(\G)$ is normal.
\end{prop}


\begin{prop}{\rm~\cite[Proposition~8]{2p3}}\label{Av}
Let $\G$ be a connected cubic edge-transitive graph and let $G\leq\Aut(\G)$ be transitive on the edges of $\G$. For any $v\in V(\G)$, the stabilizer $G_v$ has order $2^r\cdot 3$ with $r\geq 0$.
\end{prop}

\section{A few technical lemmas}

In this section, we shall give two easily proved lemmas about metacyclic $p$-groups that are useful in this paper.

\begin{lem}\label{comput}
Let $p$ be an odd prime, and let $H$ be a metacyclic $p$-group generated by $a, b$ with the following defining relations:
\begin{center}
$a^{p^m}=b^{p^n}=1, b^{-1}ab=a^{1+p^r}$,
\end{center}
where $m, n, r$ are positive integers such that $r<m\leq n+r$. Then the following hold:

\begin{enumerate}
  \item [{\rm (1)}]\ For any $i\in\mz_{p^m}, j\in\mz_{p^n}$, we have
  $$a^ib^j=b^ja^{i(1+p^r)^j}.$$

  \item [{\rm (2)}]\ For any positive integer $k$ and for any $i\in\mz_{p^m}, j\in\mz_{p^n}$, we have
  $$(b^ja^i)^k=b^{kj}a^{i\sum_{s=0}^{k-1}(1+p^r)^{sj}}.$$

  \item [{\rm (3)}]\ For any $i_1,i_2\in\mz_{p^m}$, $j_1,j_2\in\mz_{p^n}$, we have
  $$(b^{j_1}a^{i_1})(b^{j_2}a^{i_2})=b^{j_1+j_2}a^{{i_1(1+p^r)}^{j_2}+i_2}.$$
\end{enumerate}
\end{lem}
\demo
For any $i\in\mz_{p^m}, j\in\mz_{p^n}$, since $b^{-1}ab=a^{1+p^r}$, we have $b^{-j}ab^j=a^{(1+p^r)^j}$, and then $b^{-j}a^ib^j=a^{i(1+p^r)^j}$.
It follows that $a^ib^j=b^ja^{i(1+p^r)^j}$, and so (1) holds.

For any positive integer $k$ and for any $i\in\mz_{p^m}, j\in\mz_{p^n}$, if $k=1$, then $(2)$ is clearly true.
Now we assume that $k>1$ and (2) holds for any positive integer less than $k$.
Then $(b^ja^i)^{k-1}=$$b^{(k-1)j}a^{i\sum_{s=0}^{k-2}(1+p^r)^{sj}}$, and
then
\[
\begin{array}{ll}
\mbox{$(b^ja^i)^k$} & \mbox{$=b^ja^i(b^ja^i)^{k-1}$}\\
 & \mbox{$=b^ja^i[b^{(k-1)j}a^{i\sum_{s=0}^{k-2}(1+p^r)^{sj}}]$}\\
 & \mbox{$=b^j[a^ib^{(k-1)j}]a^{i\sum_{s=0}^{k-2}(1+p^r)^{sj}}$}\\
 & \mbox{$=b^j[b^{(k-1)j}a^{i(1+p^r)^{(k-1)j}}]a^{i\sum_{s=0}^{k-2}(1+p^r)^{sj}}$}\\
 & \mbox{$=b^{kj}a^{i\sum_{s=0}^{k-1}(1+p^r)^{sj}}$.}
\end{array}
\]
By induction, we have (2) holds.

For any $i_1,i_2\in\mz_{p^m}$ and $j_1,j_2\in\mz_{p^n}$, from $(1)$ it follows that $$(b^{j_1}a^{i_1})(b^{j_2}a^{i_2})=b^{j_1}(a^{i_1}b^{j_2})a^{i_2}=b^{j_1}(b^{j_2}a^{i_1(1+p^r)^{j_2}})a^{i_2}=b^{j_1+j_2}a^{{i_1(1+p^r)}^{j_2}+i_2},$$
and so (3) holds.
\hfill\qed

\begin{lem}\label{comput2}
Let $p$ be an odd prime, and let $H$ be an inner-abelian metacyclic $p$-group generated by $a, b$ with the following defining relations:
\begin{center}
$a^{p^m}=b^{p^n}=1, b^{-1}ab=a^{1+p^r}$,
\end{center}
where $m, n, r$ are positive integers such that $m\geq 2, n\geq 1$ and $r=m-1$. Then the following hold:

\begin{enumerate}
  \item [{\rm (1)}]\ For any positive integer $k$, we have
  $$a^{(1+p^r)^k}=a^{1+kp^r}.$$

  \item [{\rm (2)}]\ For any $i\in\mz_{p^m}, j\in\mz_{p^n}$, we have
  $$(b^ja^i)^p=b^{jp}a^{ip}.$$

  \item [{\rm (3)}]\ $H'\cong\z_p$.
\end{enumerate}
\end{lem}
\demo For (1), the result is clearly true if $k=1$. In what follows, assume $k\geq 2$. Since $r=m-1$ and $m\geq 2$, we have $2r\geq m$. This implies that $a^{p^{2r}}=1$, and hence $a^{p^{\ell r}}=1$ for any $\ell\geq2$. It then follows that
\[
\begin{array}{ll}
\mbox{$a^{(1+p^r)^k}$} & \mbox{$=a^{[C_k^0\cdot 1^k\cdot (p^r)^0+C_k^1\cdot 1^{k-1}\cdot (p^r)^1+C_k^2\cdot 1^{k-2}\cdot (p^r)^2+...C_k^k\cdot 1^{0}\cdot (p^r)^k]}$}\\
 & \mbox{$=a^{C_k^0\cdot (p^r)^0}\cdot a^{C_k^1\cdot (p^r)^1}\cdot a^{C_k^2\cdot (p^r)^2}\cdot ...\cdot a^{C_k^k\cdot (p^r)^k}$}\\
 & \mbox{$=a\cdot (a^{p^{r}})^{C_k^1}\cdot (a^{p^{2r}})^{C_k^2}\cdot ...\cdot (a^{p^{kr}})^{C_k^k}$}\\
 & \mbox{$=a\cdot a^{kp^r}$}\\
 & \mbox{$=a^{1+kp^r}$,}
\end{array}
\]
and so (1) holds. (Here for any integers $N\geq l\geq 0$, we denote by $C_N^l$ the binomial coefficient, that is, $C_N^l=\frac{N!}{l!(N-l)!}$.)

For (2), for any positive integer $k$ and for any $i\in\mz_{p^m}, j\in\mz_{p^n}$, by Lemma~\ref{comput}~$(2)$, we have
\[
\begin{array}{ll}
\mbox{$(b^ja^i)^{p}$} & \mbox{$=b^{jp}a^{i[1+(1+p^r)^{j}+(1+p^r)^{2j}+...+(1+p^r)^{(p-1)j}]}$}\\
 & \mbox{$=b^{jp}a^{i[1+(1+j\cdot p^r)+(1+2j\cdot p^r)+...+(1+(p-1)\cdot jp^r)]}$}\\
 & \mbox{$=b^{jp}a^{i(p+\frac{1}{2}p(p-1)\cdot jp^r)}$}\\
 & \mbox{$=b^{jp}a^{ip}$. }
\end{array}
\]
Hence $(2)$ holds.

From \cite{p-gp} we can obtain $(3)$.
\hfill\qed

\section{Proof of Theorem~\ref{5no}}

We shall prove Theorem~\ref{5no} by a series of lemmas. We first prove three lemmas regarding cubic edge-transitive graphs of order twice a prime power.

\begin{lem}\label{N}
Let $\G$ be a connected cubic edge-transitive graph of order $2p^n$ with $p$ an odd prime and $n\geq 2$.
Let $G\leq\Aut(\G)$ be transitive on the edges of $\G$.
Then any minimal normal subgroup of $G$ is an elementary abelian $p$-group.
\end{lem}

\f\demo Let $N$ be a minimal normal subgroup of $G$. If $G$ is transitive on the arcs of $\G$, then by {\rm~\cite[Lemma~3.1]{2pn}}, $N$ is an
elementary abelian $p$-group, as required.

In what follows, assume that $G$ is not transitive on the arcs of $\G$. Then since $\G$ has valency $3$, $\G$ is semisymmetric and so it is bipartite. Let $B_0$ and $B_1$ be the two partition sets of $V(\G)$. Then $B_0, B_1$ are just the two orbits $G$ on $V(\G)$ and have size $p^n$.   Recalling that $N\unlhd G$, each orbit of $N$ has size dividing $p^n$. So, if $N$ is solvable, then $N$ must be an elementary abelian $p$-group, as required.

Suppose that $N$ is non-solvable. By Proposition~\ref{Av}, we have $|G|=2^r\cdot 3\cdot p^n$, where $r\geq 0$. If $p=3$, then by Burnside $p^aq^b$-theorem, $G$ would be solvable, which is impossible because $N$ is non-solvable. Thus, $p>3$. Since $N$ is a minimal normal subgroup of $G$, $N$ is a product of some isomorphic non-abelian simple groups. Observing that $3^2\nmid |G|$, by \cite[pp.12-14]{Gorenstein}, we obtain that $N\cong A_5$ or $\PSL(2,7)$. Then $p=5$ or $7$, and $p^2\nmid |N|$. Since $n\geq 2$, it follows that $N$ is intransitive on each bipartition sets of $\G$. By Proposition~\ref{intransitive}, $N$ is semiregular on $V(\G)$, and so $|N|\ |\ p^n$, which is impossible. This completes the proof of our lemma.\hfill\qed

\begin{lem}\label{11}
Let $p\geq 5$ be a prime and let $\G$ be a connected cubic edge-transitive graph of order $2p^n$ with $n\geq 1$.
Let $A=\Aut(\G)$ and let $H$ be a Sylow $p$-subgroup of $A$. Then $\G$ is a bi-Cayley graph over $H$, and moreover,
if $p\geq 11$, then $\G$ is a normal bi-Cayley graph over $H$.
\end{lem}

\f\demo By Proposition~\ref{Av}, the stabilizer of any $v\in V(\G)$ in $A$ has order dividing $2^r\cdot 3$ with $r\geq 0$. Recalling $H$ is a Sylow $p$-subgroup of $A$, $H$ must be semiregular on $V(\G)$ since $p\geq 5$. Since $\G$ is edge-transitive, $\G$ is either arc-transitive or semisymmetric, and so $p^n\ |\ |A|$. It follows that $p^n\ |\ |H|$, and so $|H|=p^n$. Thus, $H$ has two orbits on $V(\G)$, and hence $\G$ is bi-Cayley graph over $H$.

Now suppose that $p\geq 11$. We shall prove the second assertion. It suffices to prove that $H\unlhd A$. Use induction on $n$. If $n=1$, then $\G$ is symmetric by {\cite[Theorem~2]{2p}}, and then by~\cite[Theorem~1]{BiCayley.7} (see also \cite[Table~1]{Cheng87} or \cite[Proposition~2.8]{2pn}), we have $H\unlhd A$, as required.
Assume $n\geq 2$. Take $N$ to be a minimal normal subgroup of $A$. By Lemma~\ref{N}, $N$ is an elementary abelian $p$-group and $|N|\mid p^n$.
Consider the quotient graph  $\G_N$ of $\G$ corresponding to the orbits of $N$.  If $|N|=p^n$, then $H=N\unlhd A$, as required.
Suppose that $|N|<p^n$. Then each orbit of $N$ has size at most $p^{n-1}$, and by Propositions~\ref{intransitive} and \ref{3orbits}, $N$ is semiregular, and $\G_N$ is of valency $3$ with $A/N$ as an edge-transitive group of automorphisms of $\G_N$. Clearly, $\G_N$ has order  $2p^m$ with $m<n$. By induction, we have any Sylow $p$-subgroup of $\Aut(\G_N)$ is normal. It follows that $H/N\unlhd A/N$ because $H/N$ is a Sylow $p$-subgroup of $A/N$. Therefore, $H\unlhd A$, as required.
\hfill\qed

\begin{lem}\label{Q}
Let $\G$ be a connected cubic edge-transitive graph of order $2p^n$ with $p=5$ or $7$ and $n\geq 2$.
Let $Q=O_p(A)$ be the maximal normal $p$-subgroup of $A=\Aut(\G)$. Then $|Q|=p^n$ or $p^{n-1}$.
\end{lem}
\demo Let $|Q|=p^m$ with $m\leq n$. Suppose that $n-m\geq 2$. Then by Propositions~\ref{3orbits} and \ref{intransitive},
the quotient graph $\G_Q$ is a connected cubic graph of order $2p^{n-m}$ with $A/Q$ as an edge-transitive group of automorphisms.
Let $N/Q$ be a minimal normal subgroup of $A/Q$. By Lemma~\ref{N}, $N/Q$ is an elementary abelian $p$-group.
It follows that $N\unlhd A$ and $Q<N$, contrary to the maximality of $Q$.
Thus $n-m\leq 1$, and so $|Q|=p^n$ or $p^{n-1}$.\hfill\qed


Now we are ready to consider cubic edge-transitive bi-Cayley graphs over a metacyclic $p$-group. We first prove that $p=3$.

\begin{lem}\label{57}
Let $\G$ be a connected cubic edge-transitive bi-Cayley graph over a non-abelian metacyclic $p$-group $H$ with $p$ an odd prime. Then $p=3$.
\end{lem}

\f\demo Suppose to the contrary that $p>3$. Let $A=\Aut(\G)$. Then $R(H)$ is a Sylow $p$-subgroup of $A$. We shall first prove the following claim.\medskip

\f{\bf Claim}.\ $R(H)\unlhd A$.\medskip

Suppose to the contrary that $R(H)$ is not normal in $A$. By Lemma~\ref{11}, we have $p=5$ or $7$.
Let $N$ be the maximal normal $p$-subgroup of $A$. Then $N\leq R(H)$, and by Lemma~\ref{Q}, we have $|R(H): N|=p$. Then the quotient graph $\G_N$ is a cubic graph of order $2p$ with $A/N$ as an edge-transitive automorphism group. By \cite{Sym768, Semisym768}, if $p=5$, then $\G_N$ is the Petersen graph, and if $p=7$, then $\G_N$ is the Heawood graph. Since $A/N$ is transitive on the edges of $\G_N$ and $R(H)/N$ is non-normal in $A/N$, it follows that
$$
\begin{array}{ll}
A_5\lesssim A/N\lesssim S_5, & {\rm if }\ p=5;\\
\PSL(2,7)\lesssim A/N\lesssim \PGL(2,7), & {\rm if }\ p=7.
\end{array}
$$
Let $B/N$ be the socle of $A/N$. Then $B/N$ is also edge-transitive on $\G_N$, and so $B$ is also edge-transitive on $\G$.
Let $C=C_B(N)$. By Proposition~\ref{NC}, $B/C\lesssim\Aut(N)$. And $C/(C\cap N)\cong CN/N \unlhd B/N$. Since $B/N$ is non-abelian simple, one has $CN/N=1$ or $B/N$.

Suppose first that $CN/N=1$. Then $C\leq N$, and so $C=C\cap N=C_N(N)=Z(N)$. Then $B/Z(N)=B/C\lesssim\Aut(N)$.
Since $R(H)$ is a metacyclic $p$-group, $N$ is also a metacyclic $p$-group. If $N$ is non-abelian, then by Proposition~\ref{metap} and \cite[Lemma~2.6]{pcomplement}, $\Aut(N)$ is solvable. It follows that $B/Z(N)$ is solvable, and so $B$ is solvable. This is contrary to the fact that $B/N$ is non-abelian simple.

If $N$ is abelian, then $C=Z(N)=N$.
Let $$\Aut^{\Phi}(N)=\langle\a\in\Aut(N)\mid g^{\a}\Phi(N)=g\Phi(N),\forall g\in N\rangle,$$ where $\Phi(N)$ is the Frattini subgroup of $N$.
Recall that $\Aut^{\Phi}(N)$ is a normal $p$-subgroup of $\Aut(N)$ and $\Aut(N)/\Aut^{\Phi}(N)\leq \Aut(N/\Phi(N))$ (see{~\cite{commutator}}).
Let $K/C=(B/C)\cap \Aut^{\Phi}(N)$.
Then $K/C\unlhd B/C$, and so $K\unlhd B$. It follows that
$$B/K\cong(B/C)/(K/C)\cong((B/C)\cdot \Aut^{\Phi}(N))/\Aut^{\Phi}(N)\leq\Aut(N/\Phi(N)).$$
Clearly, $K/C$ is a $p$-group. Since $C=N$, $K$ is also a $p$-group.
As $N$ is the maximal normal $p$-subgroup of $A$, $N$ is also the maximal normal $p$-subgroup of $B$. This implies that $K=N$.
If $N$ is cyclic, then $N/\Phi(N)\cong \mathbb{Z}_p$, and so $B/N=B/K\lesssim\Aut(N/\Phi(N))\cong\mathbb{Z}_{p-1}$, again contrary to the fact that $B/N$ is a non-abelian simple group. If $N$ is not cyclic, then $N/\Phi(N)\cong \mathbb{Z}_p\times\mathbb{Z}_p$.
It follows that $B/N=B/K\lesssim\Aut(N/\Phi(N))\cong\GL(2, p)$. This forces that either $A_5\leq \GL(2,5)$ with $p=5$, or $\PSL(2,7)\leq \GL(2,7)$ with $p=7$. However, each of these can not happen by Magma{~\cite{Magma}}, a contradiction.

Suppose now that $CN/N=B/N$. Since $C\cap N=Z(N)$, we have $1<C\cap N\leq Z(C)$.
Clearly, $Z(C)/(C\cap N) \unlhd C/(C\cap N)\cong CN/N$. Since $CN/N=B/N$ is non-abelian simple, $Z(C)/C\cap N$ must be trivial.
Thus $C\cap N=Z(C)$, and hence $B/N=CN/N\cong C/C\cap N=C/Z(C)$. If $C=C'$, then $Z(C)$ is a subgroup of the Schur multiplier of $B/N$. However, the Schur multiplier of $A_5$ or $\PSL(2,7)$ is $\mz_2$, a contradiction. Thus, $C\neq C'$.
Since $C/Z(C)$ is non-abelian simple, one has $C/Z(C)=(C/Z(C))'=C'Z(C)/Z(C)\cong C'/(C'\cap Z(C))$, and then we have $C=C'Z(C)$. It follows that $C''=C'$. Clearly, $C'\cap Z(C)\leq Z(C')$, and $Z(C')/(C'\cap Z(C)) \unlhd C'/(C'\cap Z(C))$. Since $C'/(C'\cap Z(C))\cong C/Z(C)$ and since $C/Z(C)$ is non-abelian simple, it follows that $Z(C')/(C'\cap Z(C))$ is trivial, and so $Z(C')=C'\cap Z(C)$.
As $C/(C\cap N)\cong CN/N$ is non-abelian, we have $C/(C\cap N)=(C/(C\cap N))'=(C/Z(C))'\cong C'/(C'\cap Z(C))=C'/Z(C')$.
Since $C'=C''$, $Z(C')$ is a subgroup of the Schur multiplier of $CN/N$. However, the Schur multiplier of $A_5$ or $\PSL(2,7)$ is $\mz_2$, forcing that $Z(C')\cong\mz_2$. This is impossible because $Z(C')=C'\cap Z(C)\leq C\cap N$ is a $p$-subgroup. Claim is proved.\medskip

If $H$ is non-split, then by Proposition~\ref{pcom}, $A$ has a normal $p$-complement $Q$. By Propositions~\ref{3orbits} and \ref{intransitive}, the quotient graph $\G_Q$ would be cubic graph of odd order, a contradiction.

Thus, $H$ is split. Then we may assume that $$H=\lg a,b\ |\ a^{p^{m}}=b^{p^{n}}=1, a^b=a^{1+p^{r}}\rg,$$ where $m, n, r$ are positive integers such that $r<m\leq m+n$.

By Claim, $R(H)\unlhd A$. Since $\G$ is edge-transitive, we assume that $\G=\Bb$. By Proposition~\ref{bicayley}, we may assume that $S=\{1, g, h\}$ with $g, h\in H$. By Proposition~\ref{bicayleyaut}, there exists $\sigma_{\a, x}\in\Aut(\G)_{1_0}$, where $\a\in\Aut(H)$ and $x\in H$, such that $\s_{\a, x}$ cyclically permutates the three elements in $\G(1_0)=\{1_1, g_1, h_1\}$. Without loss of generality, assume that $(\sigma_{\a, x})_{|\G(1_0)}=(1_1\ g_1\ h_1)$. Then $g_1=(1_1)^{\sigma_{\a, x}}=x_1$, implying that $x=g$.
Furthermore, $h_1=(g_1)^{\sigma_{\a,x}}=(gg^{\a})_1$ and $1_1=(h_1)^{\sigma_{\a, x}}=(gh^{\a})_1$.
It follows that $g^{\a}=g^{-1}h$, $h^{\a}=g^{-1}$. This implies that $\a$ is an automorphism of $H$ order dividing $3$. If $\a$ is trivial, then $h=g^{-1}$ and $g=g^{-1}h=g^{-2}$, and then $g^3=1$. Since $p>3$, we must have $h=g=1$, a contradiction. Thus, $\a$ has order $3$.  By Proposition~\ref{metap}, we must have $3\ |\ p-1$. Furthermore, $\a$ is conjugate to the following automorphism of $H$ induced by the following map:
$$\b: a\mapsto a^s, b\mapsto b,$$
where $s$ is an element of order $3$ of $\mz_{p^m}^*$.

Assume that $\b=\pi^{-1}\a\pi$ for $\pi\in\Aut(H)$. Consider the graph $\G^\pi=\BiCay(H, \emptyset, \emptyset, S^\pi)$. By Proposition~\ref{bicayleyaut}~(3), we have $\G^\pi\cong \G$, and $\s_{\b, g^{\pi}}$ cyclically permutates the three elements in $\G^\pi(1_0)=\{1_1^{\pi}, g_1^{\pi}, h_1^{\pi}\}$. For convenience of the statement, we may assume that $\pi$ is trivial and $\a=\b$.

Let $g=b^ja^i$, where $i\in\z_{p^m}$, $j\in\z_{p^n}$. Then $h=gg^{\a}=b^ja^ib^ja^{is}$.
Since $\G$ is connected, we have $H=\lg S\rg=\lg g,h\rg=$$\lg b^ja^i, b^ja^ib^ja^{is}\rg$$=\lg b^j, a^i, a^{is}\rg=\lg a^i, b^j\rg$, implying that
$i,j$ are coprime to $p$. Then there exists an integer $u$ such that $ui\equiv1\ (\mod p^m)$.
It is easy to check that the map $\g: a\mapsto a^u, b\mapsto b$ can induce an automorphism of $H$, and then $(a^i)^{\g}=a^{ui}=a$.
Again, by Proposition\ref{bicayleyaut}~(3), we have $\G\cong {\rm BiCay}(H, \emptyset, \emptyset, S^{\g})$,
where $S^{\g}=\{1, b^ja, b^jab^ja^{s}\}$. Let $\G'={\rm BiCay}(H, \emptyset, \emptyset, S^{\g})$. Then $\s_{\g^{-1}\a\g, g^\g}\in\Aut(\G')$ cyclically permutates the elements in $\G'(1_0)=\{1_1, (b^ja)_1, (b^jab^ja^s)_1\}$.

It is easy to check that $a^{\g^{-1}\a\g}=(a^i)^{\a\g}=(a^{is})^{\g}=a^s$ and $b^{\a^{\g}}=b$.
It then follows that $1_1^{\s_{\a^{\g}, b^ja}}=(b^ja)_1$, $(b^ja)_1^{\s_{\a^{\g}, b^ja}}=(b^jab^ja^s)_1$,
and $(b^jab^ja^s)_1^{\s_{\a^{\g},b^ja}}$$=(b^ja(b^jab^ja^s)^{\a^{\g}})_1$$=(b^jab^ja^sb^ja^{s^2})_1=(b^{3j}a^{(1+p^r)^{2j}+s(1+p^r)^{j}+s^2})_1\neq 1_1$. This is a contradiction. Thus $p=3$. \hfill\qed

In what follows, we consider cubic edge-transitive bi-Cayley graph over the group $H$, where $H$ is a non-abelian metacyclic $3$-group.

\begin{lem}\label{HNormal}
Let $\G=\B$ be a connected cubic edge-transitive bi-Cayley graph over a non-abelian metacyclic $3$-group $H$ with $|H|=3^s$, where $s\geq 4$.
Then $\G$ is a normal bi-Cayley graph over $H$.
\end{lem}

\f\demo Let $A=\Aut(\G)$ and let $P$ be a Sylow $3$-subgroup of $A$ such that $R(H)\leq P$. By Proposition~\ref{Av}, we have $|A|=3^{s+1}\cdot 2^r$ with $r\geq 0$. This implies that $|P|=3|R(H)|$, and so $|P_{1_0}|=|P_{1_1}|=3$. Thus, $P$ is transitive on the edges of $\G$. Clearly, $R(H)\unlhd P$. This implies that the two orbits $H_0, H_1$ of $R(H)$ do not contain the edges of $\G$, and so $R=L=\emptyset$.\medskip

\f{\bf Claim}\ $P\unlhd A$.

Let $M\unlhd A$ be maximal subject to that $M$ is intransitive on both $H_0$ and $H_1$. By Proposition~\ref{3orbits} and Proposition~\ref{intransitive}, $M$ is semiregular on $V(\G)$ and the quotient graph $\G_M$ of $\G$ relative to $M$ is a cubic graph with $A/M$ as an edge-transitive group of automorphisms. Assume that $|M|=3^t$. Then $|V(\G_M)|=2\cdot 3^{s-t}$. If $s-t\leq 2$, then by \cite{Sym768, Semisym768}, $\G_M$ is isomorphic to $\F006A$ or the Pappus graph $\F018A$, and then $\Aut(\G_M)$ has a normal Sylow $3$-subgroup. It follows that $P/M\unlhd A/M$, and so $P\unlhd A$, as claimed.

Now assume that $s-t>2$. Take a minimal normal subgroup $N/M$ of $A/M$. By Lemma~\ref{N}, $N/M$ is an elementary abelian $3$-group. By the maximality of $M$, $N$ is transitive on at least one of $H_0$ and $H_1$, and so $3^{s}\ |\ |N|$. If $3^{s+1}\ |\ |N|$, then $P=N\unlhd A$, as claimed. Assume that $|N|=3^s$. If $N$ is transitive on both $H_0$ and $H_1$, then $N$ is semiregular on both $H_0$ and $H_1$, and then $\G_M$ would be a cubic bi-Cayley graph on $N/M$. Since $\G_M$ is connected, by Proposition~\ref{bicayley}, $N/M$ is generated by two elements, and so $N/M\cong\mz_3$ or $\mz_3\times\mz_3$. This implies that $|V(\G_M)|=6$ or $18$, contrary to the assumption that $|V(\G_M)|=2\cdot 3^{s-t}>18$. Thus, we may assume that $N$ is transitive on $H_0$ but intransitive on $H_1$. Then $N/M\neq R(H)M/M$, and so $NR(H)M/M=P/M$. Since $|P/M: R(H)M/M|\ |\ 3$, one has $|N/M: (N/M\cap R(H)M/M)|\ |\ 3$, and since $H$ is metacyclic, one has $N/M\cap R(H)M/M$ is also metacyclic and so is a two-generator group. This implies that $|N/M|\ |\ 3^3$, and so $|N/M|=3^3$ because $|N/M|=3^{s-t}>9$. Then $|V(\G_M)|=2\cdot |N/M|=54$. Since $s\geq 4$, we have $|M|\geq 3$. If $M\nleq R(H)$, then $P=MR(H)$ and then $N/M\leq R(H)M/M$. As $H$ is metacyclic, $N/M$ is also metacyclic, and so $|N/M|=3$ or $9$, a contradiction. Thus, $M\leq R(H)$, and hence $M$ is metacyclic. Then $M/\Phi(M)\cong\mz_3$ or $\mz_3\times\mz_3$. Since $\Phi(M)$ is characteristic in $M$, one has $\Phi(M)\unlhd A$ because $M\unlhd A$. Then the quotient graph $\G_{\Phi(M)}$ is a cubic graph of order $2\cdot 3^4$ or $2\cdot 3^5$ with $A/\Phi(M)$ as an edge-transitive group of automorphisms. By \cite{Sym768,Semisym768} and Magma~\cite{Magma}, we obtain that every Sylow $3$-subgroup of $\Aut(\G_{\Phi (M)})$ is normal. This implies that $P/\Phi(M)\unlhd A/\Phi(M)$, and so $P\unlhd A$, completing the proof of our claim.\medskip

Now we are ready to finish the proof of our lemma. By Claim, we have $P\unlhd A$. Since $|P: R(H)|=3$, one has $\Phi(P)\leq R(H)$. As $H$ is non-abelian, one has $\Phi (P)<R(H)$ for otherwise, we would have $P$ is cyclic and so $H$ is cyclic which is impossible. Then $\Phi(P)$ is intransitive on both $H_0$ and $H_1$, the two orbits of $R(H)$ on $V(\G)$. Since $\Phi(P)$ is characteristic in $P$, $P\unlhd A$ gives that $\Phi(P)\unlhd A$. By Propositions~\ref{3orbits} and \ref{intransitive}, the quotient graph $\G_{\Phi(P)}$ of $\G$ relative to $\Phi(P)$ is a cubic graph with $A/\Phi(P)$ an edge-transitive group of automorphisms. Furthermore, $P/\Phi(P)$ is transitive on the edges of $\G_{\Phi(P)}$. Since $P/\Phi(P)$ is abelian, it is easy to see that $\G_{\Phi(P)}\cong K_{3,3}$, and so $P/\Phi(P)\cong\mz_3\times\mz_3$. Since $|P|=3^{s+1}\geq 3^5$, one has $|\Phi(P)|=3^{s-1}\geq 3^3$.

Let $\Phi_2$ be the Frattini subgroup of $\Phi(P)$.
Then $\Phi_2\unlhd A$ because $\Phi_2$ is characteristic in $\Phi (P)$ and $\Phi(P)\unlhd A$.
Clearly, $\Phi_2\leq \Phi(P)<R(H)$, so $\Phi_2$ is intransitive on both $H_0$ and $H_1$.
Consider the quotient graph $\G_{\Phi_2}$ of $\G$ relative to $\Phi_2$. By Propositions~\ref{3orbits} and \ref{intransitive}, $\G_{\Phi_2}$ is a cubic graph with $A/\Phi_2$ as an edge-transitive group of automorphisms. Furthermore, $\G_{\Phi_2}$ is a bi-Cayley graph over the group $R(H)/\Phi_2$. Again, since $H$ is a metacyclic group, we have $\Phi(P)/\Phi_2\cong\mz_3$ or $\mz_3\times\mz_3$. If $\Phi(P)/\Phi_2\cong\mz_3$, then $\Phi(P)$ is a cyclic $3$-group, and so $\G$ is an edge-transitive cyclic cover of $\G_{\Phi(P)}\cong K_{3,3}$. By Feng et al.~\cite{K33Feng, K33Wang}, we have $\G$ is isomorphic to either $K_{3,3}$ or the Pappus graph, a contradiction.

Thus, $\Phi(P)/\Phi_2\cong\mz_3\times\mz_3$. Since $|\Phi(P)|=3^{s-1}\geq 3^3$, one has $|\Phi_2|\geq 3$.
Let $\Phi_3$ be the Frattini subgroup of $\Phi_2$. Then $\Phi_3$ is characteristic in $\Phi_2$, and so normal in $A$ because $\Phi_2\unlhd A$.
As $\Phi_2\leq R(H)$, one has $\Phi_2/\Phi_3\cong\z_3$ or $\z_3\times\z_3$, and so $|R(H)/\Phi_3|=3^4$ or $3^5$. Clearly, $\Phi_3$ is intransitive on both $H_0$ and $H_1$. Again, by Propositions~\ref{3orbits} and \ref{intransitive}, the quotient graph $\G_{\Phi_3}$ is a cubic graph of order $162$ or $486$ with $A/\Phi_3$ as an edge-transitive group of automorphisms. Observe that $R(H)/\Phi_3$ is metacyclic semiregular on $V(\G_{\Phi_3})$ with two orbits.

If $|\G_{\Phi_3}|=486$, then by \cite{Sym768, Semisym768}, $\G_{\Phi_3}$ is semisymmetric or symmetric. For the former, by Magma \cite{Magma}, all semiregular subgroups of $\Aut(\G_{\Phi_2})$ of order $243$ are normal, and so $R(H)/\Phi_3\unlhd \Aut(\G_{\Phi_3})$. It follows that $R(H)/\Phi_3\unlhd A/\Phi_3$, and so $R(H)\unlhd A$, as required. If $\G_{\Phi_3}$ is symmetric, then by \cite{Sym768}, $\G_{\Phi_3}\cong \F486A$, $\F486B$, $\F486C$ or $\F486D$. By Magma \cite{Magma}, if $\G_{\Phi_3}\cong \F486B, \F486C$ or $\F486D$, then $\Aut(\G_{\Phi_3})$ does not have a metacyclic semiregular subgroup of order $243$, a contradiction.
If $\G_{\Phi_3}\cong \F486A$, then by Magma \cite{Magma}, all semiregular subgroups of $\Aut(\G_{\Phi_3})$ of order $243$ are normal, and so $R(H)/\Phi_3\unlhd \Aut(\G_{\Phi_3})$. It follows that $R(H)/\Phi_3\unlhd A/\Phi_3$, and so $R(H)\unlhd A$, as required.

If $|\G_{\Phi_3}|=162$, then by \cite{Sym768, Semisym768}, $\G_{\Phi_3}$ is symmetric, and is isomorphic to $\F162A$, $\F162B$ or $\F162C$. By Magma \cite{Magma}, if $\G_{\Phi_3}\cong \F162C$, then $\Aut(\G_{\Phi_3})$ does not have a metacyclic semiregular subgroup of order $81$, a contradiction. If $\G_{\Phi_3}\cong \F162A$ or $\F162B$, then by Magma \cite{Magma}, all semiregular subgroups of $\Aut(\G_{\Phi_3})$ of order $81$ are normal, and so $R(H)/\Phi_3\unlhd \Aut(\G_{\Phi_3})$. It follows that $R(H)/\Phi_3\unlhd A/\Phi_3$, and so $R(H)\unlhd A$, as required.
\hfill\qed

\f{\bf Proof of Theorem~\ref{5no}}\ Let $\G=\B$ be a connected cubic edge-transitive bi-Cayley graph over a non-abelian metacyclic $p$-group $H$ with $p$ an odd prime. By Lemma~\ref{57}, we have $p=3$, and since $H$ is a non-abelian metacyclic $3$-group, we have $|H|=3^s$ with $s\geq 3$. If $s=3$, then $\G$ has order $54$, and by \cite{Sym768, Semisym768}, $\G$ is isomorphic to $\F054$ or the Gray graph. However, by Magma \cite{Magma}, $\Aut(\F054)$ does not have a non-abelian metacyclic $3$-subgroup which acts semiregularly on the vertex set of $\F054$ with tow orbits. It follows that $\G$ is isomorphic to Gray graph. If $s>3$, then by Lemma~\ref{HNormal}, $R(H)\unlhd\Aut(\G)$, as required.\hfill\qed

\section{A class of cubic edge-transitive bi-$3$-metacirculants}

In this section, we shall use Theorem~\ref{5no} to give a characterization of connected cubic edge-transitive bi-Cayley graphs over inner-abelian metacyclic $3$-groups.

\subsection{Construction}\label{sec-5}

We shall first construct two classes of connected cubic edge-transitive bi-Cayley graphs over inner-abelian metacyclic $3$-groups.

\begin{construction}\label{con-1}
Let $t$ be a positive integer, and let $$\GG=\lg a, b\ |\ a^{3^{t+1}}=b^{3^t}=1, b^{-1}ab=a^{1+3^t}\rg.$$
Let $S=\{1, a, a^{-1}b\}$, and let $\G_{t}=\BiCay(\GG, \emptyset, \emptyset, S)$.
\end{construction}

\begin{lem}\label{Example1}
For any integer $t$, the graph ${\G_t}$ is semisymmetric.
\end{lem}

\f\demo We first prove the following four claims.\medskip

\f{\bf Claim 1.} $\GG$ has an automorphism $\a$ mapping $a, b$ to $a^{-2}b, a^{3^t-3}b$, respectively.\medskip

Let $x=a^{-2}b$ and $y=a^{3^t-3}b$. Then,
$$\begin{array}{l}
(yx^{-1})^{3^t+1}=[(a^{3^t-3}b)(a^{-2}b)^{-1}]^{3^t+1}=(a^{3^t-1})^{3^t+1}=a^{-1},\\
((yx^{-1})^{3^t+1})^{-2}\cdot x=a^2\cdot a^{-2}b=b,
\end{array}
$$
and hence $\lg a, b\rg=\lg x, y\rg$.

By Lemma~\ref{comput2}~(2), we have $x^{3^{t+1}}=(a^{-2}b)^{3^{t+1}}=1$ and $y^{3^{t}}=(a^{3^t-3}b)^{3^{t}}=1$.
Furthermore, we have
$$x^{1+3^t}=(a^{-2}b)^{1+3^t}=(a^{-2}b)(a^{-2}b)^{3^t}=a^{-2}ba^{-2\cdot 3^t}=a^{-2-2\cdot 3^t}b=a^{3^t-2}b,$$
and
$$
\begin{array}{lll}
y^{-1}xy&=&(a^{3^t-3}b)^{-1}(a^{-2}b)(a^{3^t-3}b)\\
&=&(b^{-1}a^{3-3^t}a^{-2}b)a^{3^t-3}b\\
&=&(b^{-1}a^{1-3^t}b)a^{3^t-3}b\\
&=&a^{(1+3^t)(1-3^t)}a^{3^t-3}b\\
&=&a^{3^t-2}b\\
&=&x^{1+3^t}.
\end{array}$$
It follows that $x$ and $y$ have the same relations as do $a$ and $b$.
Thus, the map $\a: a\mapsto a^{-2}b, b\mapsto a^{3^t-3}b$ induces an automorphism of $\GG$, as claimed.\medskip

\f{\bf Claim 2.}\ $\GG$ has no automorphism mapping $a, b$ to $a^{-1}, a^{3t}b^{-1}$, respectively.\medskip

Suppose to the contrary that $\GG$ has an automorphism, say $\b$, such that $a^\b=a^{-1}, b^\b=a^{3t}b^{-1}$.
Then $(b^{-1}ab)^\b=(a^{3^t+1})^\b$, and so
$$\begin{array}{lll}
a^{-3^t-1}&=&(a^{3^t+1})^\b=(b^{-1}ab)^\b\\
&=&(a^{3^t}b^{-1})^{-1}\cdot a^{-1}\cdot (a^{3^t}b^{-1})\\
&=&ba^{-1}b^{-1}=a^{-(1+3^t)^{3^t-1}}=a^{-1+3^t}.
\end{array}
$$
It follows that $a^{2\cdot 3^t}=1$, and so $3^{t+1}\ |\ 2\cdot 3^t$, a contradiction.\medskip

\f{\bf Claim 3.}\ $\GG$ has no automorphism mapping $a, b$ to $b^{-1}a, b^{-1}$, respectively.\medskip

Suppose to the contrary that there exists $\g\in\Aut(\GG)$ such that $a^\g=b^{-1}a, b^\g=b^{-1}$.
Then $(b^{-1}ab)^\g=(a^{1+3^t})^\g$, and then
$$
\begin{array}{lll}
b^{-1}a^{3^t+1}=(b^{-1}a)^{1+3^t}=(a^{1+3^t})^\g=(b^{-1}ab)^\g=b(b^{-1}a)b^{-1}=ab^{-1}.\\
\end{array}
$$
It follows that $b^{-1}a^{3^t+1}b=a$, and so $a^{3^{2t}+2\cdot 3^t+1}=a^{2\cdot 3^t+1}=a$, forcing that $3^{t+1}\ |\ 2\cdot 3^t$, a contradiction.\medskip

Now we are ready to finish the proof. By Claim~1, there exists $\a\in\Aut(\GG)$ such that $a^{\a}=a^{-2}b$ and $b^{\a}=a^{3^t-3}b$.
Then $(a^{-1}b)^{\a}=(a^{-2}b)^{-1}(a^{3^t-3}b)=b^{-1}a^{3^t-1}b=a^{-1}$. It then follows that
$$S^{\a}=\{1^{\a},a^{\a},(a^{-1}b)^{\a}\}=\{1,a^{-2}b,a^{-1}\}=a^{-1}S.$$ By Proposition~\ref{bicayleyaut},
$\s_{\a, a}$ is an automorphism of $\G_t$ fixing $1_0$ and cyclically permutating the three neighbors of $1_0$.
Set $B=R(\GG)\rtimes\lg \s_{\a,a}\rg$. Then $B$ acts regularly on the edges of ${\G_t}$.

If $t=1$, then by Magma~\cite{Magma}, $\G_1$ is isomorphic to the Gray graph, which is semisymmetric.
In what follows, assume that $t>1$. By Theorem~\ref{5no}, $\G_t$ is a normal bi-Cayley graph over $R(H)$. Suppose that $\G_t$ is vertex-transitive.
Then $\G_t$ is also arc-transitive. So, there exist $f\in\Aut(\GG), g, h\in\GG$ so that $\d_{f, g, h}$ is an automorphism of $\G_t$ taking the arc $(1_0, 1_1)$ to $(1_1, 1_0)$. By the definition of $\d_{f, g, h}$, one may see that $g=h=1$ and $S^{f}=S^{-1}$, namely,
$$\{1, a, a^{-1}b\}^f=\{1, a^{-1}, b^{-1}a\}.$$
So, $f$ takes $(a, a^{-1}b)$ either to $(a^{-1}, b^{-1}a)$ or to $(b^{-1}a, a^{-1})$. However, this is impossible by Claims 2-3. Therefore, ${\G_t}$ is semisymmetric.\hfill\qed

\begin{construction}\label{con-2}
Let $t$ be a positive integer, and let
$$\H=\lg a, b\ |\ a^{3^{t+1}}=b^{3^{t+1}}=1, b^{-1}ab=a^{1+3^t}\rg.$$
Let $T=\{1,b,b^{-1}a\}$, and let $\Sigma_{t}=\BiCay(\H,\emptyset, \emptyset, T)$.
\end{construction}

\begin{lem}\label{Example2}
For any positive integer $t$, the graph $\Sigma_t$ is symmetric.
\end{lem}

\f\demo We first prove the following two claims.\medskip

\f{\bf Claim 1} $\H$ has an automorphism $\a$ mapping $a, b$ to $a^{2\cdot 3^t+1}b^{-3}, a^{2\cdot 3^t+1}b^{-2}$, respectively.\medskip

Let $x=a^{2\cdot 3^t+1}b^{-3}$ and $y=a^{2\cdot 3^t+1}b^{-2}$. Note that $((y^{-1}x)^{-1})=b$ and $xb^3=a^{2\cdot 3^t+1}$. This implies that $\lg x, y\rg=\lg a, b\rg=\H$.

By Lemma~\ref{comput2} (2), we have $x^{3^{t+1}}=(a^{-2}b)^{3^{t+1}}=1$ and $y^{3^{t+1}}=(a^{3^t-3}b)^{3^{t}}=1$. Furthermore, we have
$$
\begin{array}{lll}
y^{-1}xy&=&(a^{2\cdot 3^t+1}b^{-2})^{-1}(a^{2\cdot 3^t+1}b^{-3})(a^{2\cdot 3^t+1}b^{-2})\\
&=& b^{-1}a^{2\cdot 3^t+1}b^{-2}= b^{-1}a^{2\cdot 3^t+1}bb^{-3}\\
&=& a^{(2\cdot 3^t+1)(3^t+1)}b^{-3}=ab^{-3}=x^{3^t}x\\
&=& x^{3^t+1}.
\end{array}
$$
It follows that $x$ and $y$ have the same relations as do $a$ and $b$. Therefore, $\H$ has an automorphism taking $(a, b)$ to $(x, y)$, as claimed.\medskip

\f{\bf Claim 2.} $H_t$ has an automorphism $\b$ mapping $a,b$ to $a^{-1},a^{-1}b$.\medskip

Let $x=a^{-1}$ and $y=a^{-1}b$. Clearly, $\lg a,b\rg=\lg x,y\rg$. By Lemma~\ref{comput2} (2), we have that $x^{3^{t+1}}=(a^{-1})^{3^{t+1}}=1$ and $y^{3^{t+1}}=(a^{-1}b)^{3^{t+1}}=1$. Furthermore, we have
$$
\begin{array}{lll}
y^{-1}xy&=&(a^{-1}b)^{-1}(a^{-1})(a^{-1}b)=b^{-1}a^{-1}b=a^{-3^t-1}=x^{3^t+1}.
\end{array}
$$
It follows that $x$ and $y$ have the same relations as do $a$ and $b$. Therefore, $\H$ has an automorphism $\b$ which takes $(a, b)$ to $(a^{-1}, a^{-1}b)$, as claimed.\medskip

Now we are ready to finish the proof. By Claim~1, there exists $\a\in\Aut(\H)$ such that $a^{\a}=a^{2\cdot 3^t+1}b^{-3}$ and $b^{\a}=a^{2\cdot 3^t+1}b^{-2}$. Then
$$S^{\a}=\{1, b, b^{-1}a\}^\a=\{1, a^{2\cdot 3^t+1}b^{-2}, b^{-1}\}.$$
By an easy computation, we have $a^{2\cdot 3^t+1}b^{-2}=a^{2\cdot 3^t+1}b^{-3}b=b^{-3}a^{2\cdot 3^t+1}b=b^{-2}b^{-1}a^{2\cdot 3^t+1}b=b^{-2}a^{(2\cdot 3^t+1)(3^t+1)}=b^{-2}a.$ It follows that $$b^{-1}S=b^{-1}\{1, b, b^{-1}a\}=\{b^{-1}, 1, b^{-2}a\}=S^\a.$$
By Proposition~\ref{bicayleyaut}, $\s_{\a, b}$ is an automorphism of $\Sigma_t$ fixing $1_0$ and cyclically permutating the three neighbors of $1_0$.
Set $B=R(\H)\rtimes\lg \s_{\a,b}\rg$. Then $B$ acts transitively on the edges of $\Sigma_t$.

By Claim~2, there exists $\b\in\Aut(\H)$ such that $a^{\b}=a^{-1}$ and $b^{\b}=a^{-1}b$.
Then $S^{\b}=\{1, b, b^{-1}a\}^\b=\{1, a^{-1}b, b^{-1}\}=S^{-1}$. By Proposition~\ref{bicayleyaut}, $\d_{\b,1,1}$ is an automorphism of $\Sigma_t$ swapping $1_0$ and $1_1$. Thus, $\Sigma_t$ is vertex-transitive, and so $\Sigma_t$ is symmetric.
\hfill\qed

\subsection{Classification}

In this section, we shall give a classification of cubic edge-transitive bi-Cayley graph over an inner-abelian metacyclic $3$-group.

\begin{lem}\label{classify}
Let $H$ be an inner-abelian metacyclic $3$-group, and let $\G$ be a connected cubic edge-transitive bi-Cayley graph over $H$.
Then $\G\cong \G_t$ or $\Sigma_t$.
\end{lem}

\f\demo Since $H$ is an inner-abelian metacyclic $3$-group, it has order at least $3^3$. If $|H|=3^3$, then $|\G|=54$ and by \cite{Sym768, Semisym768}, we know that $\G$ is isomorphic to $\G_1$. In what follows, assume that $|H|>3^3$. By Theorem~\ref{5no}, $\G$ is a normal bi-Cayley graph over $H$. Let $\G=\BiCay(H, R, L, S)$. Since $\G$ is edge-transitive, the two orbits $H_0, H_1$ of $R(H)$ on $V(\G)$ do not contain edges, and so $R=L=\emptyset$. By Proposition~\ref{bicayley}, we may assume that $S=\{1, x, y\}$ for $x, y\in H$. Since $\G$ is connected, by Proposition~\ref{bicayley}, we have $H=\lg S\rg=\lg x, y\rg$.

Let $A=\Aut(\G)$, since $\G$ is normal and since $\G$ is edge-transitive, by Proposition~\ref{bicayleyaut}, there exists $\sigma_{\a, h}\in A_{1_0}$, where $\a\in\Aut(H)$ and $h\in H$, such that $\s_{\a, h}$ cyclically permutates the three elements in $\G(1_0)=\{1_1, x_1, y_1\}$. Without loss of generality, assume that $(\sigma_{\a, h})_{|\G(1_0)}=(1_1\ x_1\ y_1)$. Then $x_1=(1_1)^{\sigma_{\a, h}}=h_1$, implying that $x=h$.
Furthermore, $y_1=(x_1)^{\sigma_{\a, h}}=(xx^{\a})_1$ and $1_1=(y_1)^{\sigma_{\a, h}}=(xy^{\a})_1$.
It follows that $x^{\a}=x^{-1}y$ and $y^{\a}=x^{-1}$. This implies that $\a$ is an automorphism of $H$ order dividing $3$. If $\a$ is trivial, then $x=y^{-1}$ and $x=x^{-1}y=y^2$, and then $y^3=1$ and $x^3=1$. This implies that $H\cong\mz_3$ or $\mz_3\times\mz_3$, contrary to the assumption that $|H|>3^3$. Thus, $\a$ has order $3$.

Since $H$ is an inner-abelian $3$-group, by elementary group theory (see also \cite{inner}), we may assume that
$$H=\lg a, b\ |\ a^{3^{t+1}}=b^{3^s}=1, b^{-1}ab=a^{3^t+1}\rg,$$
where $t\geq 2, s\geq 1$. We first prove the following claim.\medskip

\f{\bf Claim 1}\ $H/H'=\lg aH'\rg\times\lg bH'\rg\cong \z_{3^{t}}\times\z_{3^{t}}$, $\z_{3^{t}}\times\z_{3^{t-1}}$ or $\z_{3^{t}}\times\z_{3^{t+1}}$.\medskip

By Lemma~\ref{comput2}~(3), we have the derived subgroup $R(H)'$ of $R(H)$ is isomorphic to $\z_3$. Since $R(H)'$ is characteristic in $R(H)$, $R(H)\unlhd A$ gives that $R(H)'\unlhd A$. Consider the quotient graph $\G_{R(H)'}$ of $\G$ relative to $R(H)'$. Clearly, $R(H)'$ is intransitive on both $H_0$ and $H_1$, the two orbits of $R(H)$ on $V(\G)$. By Propositions~\ref{3orbits} and \ref{intransitive}, $\G_{R(H)'}$ is a cubic graph with $A/R(H)'$ as an edge-transitive group of automorphisms. Clearly, $\G_{R(H)'}$ is a bi-Cayley graph over the abelian group $R(H)/R(H)'$.
Since $R(H)/R(H)'\unlhd A/R(H)'$, by Proposition~\ref{abelian-biCayley}, we have $R(H)/R(H)'\cong\z_{3^{m+n}}\times\z_{3^{m}}$ for some integers $m, n$ satisfying the equality ${\ld}^2-\ld+1\equiv 0\ ({\rm mod}\ 3^n)$ with $\ld\in\z_{3^n}^{*}$. This implies that $n=0$ or $1$, and so $R(H)/R(H)'\cong \z_{3^m}\times\z_{3^{m}}$ or $\z_{3^{m+1}}\times\z_{3^{m}}$.

Since $a^{3^t}=[a, b]$, one has $\lg aH'\rg\cong\mz_{3^t}$, and since $H'\cap\lg b\rg=1$, one has $H/H'=\lg aH'\rg\times\lg bH'\rg\cong\mz_{3^t}\times\mz_{3^s}$. So, if $R(H)/R(H)'\cong \z_{3^m}\times\z_{3^{m}}$, then we have $m=s=t$ , and if $R(H)/R(H)'\cong \z_{3^{m+1}}\times\z_{3^{m}}$, then $(t, s)=(m, m+1)$ or $(m+1, m)$. Claim~1 is proved. \medskip

For any $h\in H$, denote by $o(h)$ the order of $h$. Let $n={\rm Max}\{t+1, s\}$. By Lemma~\ref{comput2}~(2), it is easy to see that $3^n$ is the exponent of $H$. \medskip

\f{\bf Claim 2}\ $o(x)=o(y)=o(x^{-1}y)=3^{n}$ and $x^{3^{n-1}}\neq y^{3^{n-1}}$. \medskip

Observing that $x^{\a}=x^{-1}y$ and $y^{\a}=x^{-1}$, we have $o(x)=o(y)=o(x^{-1}y)$. By Lemma~\ref{comput2}~(2), we must have $o(x)=o(y)=o(x^{-1}y)=3^{n}$. Then $(x^{-1}y)^{3^{n-1}}\neq 1$, and again by Lemma~\ref{comput2}~(2), we have $x^{-3^{n-1}}y^{3^{n-1}}\neq 1$, namely,  $x^{3^{n-1}}\neq y^{3^{n-1}}$, as claimed.\medskip

By Claim~1, we shall consider the following three cases:\medskip

\f{\bf Case~1}\ $H/H'=\lg aH'\rg\times\lg bH'\rg\cong \z_{3^{t}}\times\z_{3^{t}}$.\medskip

In this case, we have $s=t$. By Claim~2, we have $o(x)=o(y)=o(x^{-1}y)=3^{t+1}$ and $x^{3^t}\neq y^{3^t}$.
As $H'\cong\mz_3$, we have $H'=\lg x^{3^t}\rg=\lg y^{3^t}\rg$, implying that $y^{3^t}=x^{-3^t}$.
Thus $(xy)^{3^t}=x^{3^t}y^{3^t}=x^{3^t}x^{-3^t}=1$.
Since $[x, y]\in H'$ and $H'=\lg x^{3^t}\rg$, we have $[x, y]=x^{3^t}$ or $x^{-3^t}$.
It follows that $(xy)^{-1}\cdot x\cdot(xy)=y^{-1}xy=x^{1+3^t}$ or $x^{1-3^t}$.

If $(xy)^{-1}\cdot x\cdot(xy)=y^{-1}xy=x^{1+3^t}$, then
$$H=\lg x, xy\ |\ x^{3^{t+1}}=(xy)^{3^t}=1, (xy)^{-1}\cdot x\cdot(xy)=x^{1+3^t}\rg, $$
and $S=\{1, x, y\}=\{1, x, x^{-1}(xy)\}$. So, $\G$ is isomorphic to $\G_t$ (see Construction~1).

If $(xy)^{-1}\cdot x\cdot(xy)=y^{-1}xy=x^{1-3^t}$, then
$$H=\lg x, (xy)^{-1}\ |\ x^{3^{t+1}}=[(xy)^{-1}]^{3^t}=1, (xy)^{}\cdot x\cdot(xy)^{-1}=x^{1+3^t}\rg, $$
and $S=\{1, x, y\}=\{1, x, x^{-1}[(xy)^{-1}]^{-1}\}$. By Proposition~\ref{bicayley}~(4), we have
$$\G=\BiCay(H, \emptyset, \emptyset, S)\cong \BiCay(H, \emptyset, \emptyset, S^{-1}).$$
Note that $S^{-1}=\{1, x^{-1}, y^{-1}\}=\{1, x^{-1}, (xy)^{-1}x\}$. It is easy to check that the map
$$f: x\mapsto x^{-1}, (xy)^{-1}\mapsto (xy)^{-1}x^{-3^t}$$
induces an automorphism of $H$ such that $\{1, x, x^{-1}(xy)^{-1}\}^f=S^{-1}$. By Proposition~\ref{bicayley}~(3), we have
$$\G\cong\BiCay(H, \emptyset, \emptyset, S^{-1})\cong\BiCay(H, \emptyset, \emptyset, \{1, x, x^{-1}(xy)^{-1}\})\cong\G_t,$$
as required.\medskip

\f{\bf Case 2}\ $H/H'=\lg aH'\rg\times\lg bH'\rg\cong \z_{3^{t}}\times\z_{3^{t-1}}$.\medskip

In this case, we have $s=t-1$. Let $T=\lg R(h)\ |\ h\in H, h^{3^{t-1}}=1\rg$. Then $T=\lg R(a)^9\rg\times\lg R(b)\rg$ and $T$ is characteristic in $R(H)$, and so normal in $A$ for $R(H)\unlhd A$. Furthermore, $R(H)/T\cong\mz_9$. By Propositions~\ref{3orbits} and \ref{intransitive}, the quotient graph $\G_T$ of $\G$ relative to $T$ is a cubic edge-transitive graph of order $18$. Clearly, $R(H)/T$ is semiregular on $V(\G_{T})$ with two orbits, so $\G_{T}$ is a bi-Cayley graph over the cyclic group $R(H)/T$ of order $9$. Since $R(H)/T\unlhd A/T$, by Proposition~\ref{abelian-biCayley}, there exists $\ld\in\z_{3^2}^{*}$ such that ${\ld}^2-\ld+1\equiv 0\ ({\rm mod}\ 3^2)$, which is impossible.\medskip

\f{\bf Case 3} $H/H'=\lg aH'\rg\times\lg bH'\rg\cong \z_{3^{t}}\times\z_{3^{t+1}}$.\medskip

In this case, we have $s=t+1$. Let $N=\lg h\ |\ h\in H, h^{3}=1\rg$. Then $N=\lg a^{3^t}, b^{3^t}\rg\cong\mz_3\times\mz_3$.
By Claim~2, we have $o(x)=o(y)=3^{t+1}$. Since $H=\lg x, y\rg$, one has $N=\lg x^{3^t}, y^{3^t}\rg$. As $H'\cong\mz_3$, one has $H'\leq N$.
So, $H'=\lg x^{3^t}\rg$, $\lg y^{3^t}\rg$, $\lg (xy)^{3^t}\rg$ or $\lg (xy^{-1})^{3^t}\rg$.

Recall that $H$ has an automorphism $\a$ taking $(x, y)$ to $(x^{-1}y, x^{-1})$. Suppose that one of the three subgroups: $\lg x\rg, \lg y\rg, \lg x^{-1}y\rg$ is normal in $H$. Then all of them are normal in $H$. So $H=\lg x, y\rg=\lg x\rg\times\lg y\rg$ because $|H|=3^{2(t+1)}$. This is impossible because $H$ is non-abelian. Thus, all of the three subgroups: $\lg x\rg, \lg y\rg, \lg x^{-1}y\rg$ are not normal in $H$.

It then follows that $H'=\lg (xy)^{3^t}\rg$. Then either $x^{-1}(xy)x=(xy)^{1+3^t}$ or $x^{-1}(xy)x=(xy)^{1-3^t}$.
For the former, we have
$$H=\lg xy, x\ |\ (xy)^{3^{t+1}}=x^{3^{t+1}}=1, x^{-1}(xy)x=(xy)^{3^t+1}\rg,$$
and $S=\{1, x, y\}=\{1, x, x^{-1}(xy)\}$. Hence, $\G\cong\Sigma_t$ (see Construction~2).

For the latter, we have
$$H=\lg xy, x^{-1}\ |\ (xy)^{3^{t+1}}=x^{-3^{t+1}}=1, x^{}(xy)x^{-1}=(xy)^{3^t+1}\rg,$$
and $S=\{1, x, y\}=\{1, (x^{-1})^{-1}, x^{-1}(xy)\}$. By Proposition~\ref{bicayley}~(4), we have
$$\G=\BiCay(H, \emptyset, \emptyset, S)\cong \BiCay(H, \emptyset, \emptyset, S^{-1}).$$
Note that $S^{-1}=\{1, x^{-1}, y^{-1}\}=\{1, x^{-1}, (xy)^{-1}x\}$. It is easy to check that the map
$$f': x^{-1}\mapsto x^{-1}, xy\mapsto (xy)^{3^t-1}$$
induces an automorphism of $H$ such that $\{1, x^{-1},  x(xy)\}^{f'}=S^{-1}$. By Proposition~\ref{bicayley}~(3), we have
 $$\G\cong\BiCay(H, \emptyset, \emptyset, S^{-1})\cong\BiCay(H, \emptyset, \emptyset, \{1, x^{-1},  x(xy)\})\cong\Sigma_t,$$
 as required.\hfill\qed

\section{Proof of Corollary~\ref{cor2p3}}

Let $p$ be a prime, and let $\G$ be a connected cubic edge-transitive graph of order $2p^3$. By \cite{2p}, the smallest semisymmetric graph has $20$ vertices. So, if $p=2$, then $\G$ is vertex-transitive. If $p=3$, then by \cite{Sym768,Semisym768}, we know that $\G$ is not vertex-transitive if and only if it isomorphic to the Gray graph.

Now assume that $p>3$. By Lemma~\ref{11}, $\G$ is a bi-Cayley graph over a group $H$ of order $p^3$.
Suppose that $\G$ is not vertex-transitive. Then $\G$ is bipartite with the two orbits of $H$ as its two parts. So we may let $\G=\BiCay(H, \emptyset, \emptyset, S)$. By Proposition~\ref{bicayley}, we may assume that $S=\{1, a, b\}$ form $a,b\in H$. If $H$ is abelian, then $H$ has an automorphism $\a$ which maps every element of $H$ to its inverse. By Proposition~\ref{bicayleyaut}, $\d_{\a,1,1}$ is an automorphism $\G$ swapping the two parts of $\G$, and so $\G$ is vertex-transitive, a contradiction. If $H$ is non-abelian, then $H$ is either metacyclic or isomorphic to the following group:
$$J=\lg a,b,c\ |\ a^p=b^p=c^p=1, c=[a,b], [a,c]=[b,c]=1\rg.$$
By Theorem~\ref{5no}, $H$ is non-metacyclic. If $H\cong J$, then it is easy to see that $J$ has an automorphism taking $(a, b)$ to $(a^{-1}, b^{-1})$. Again, by Proposition~\ref{bicayleyaut}, $\d_{\a,1,1}$ is an automorphism of $\G$ swapping the two parts of $\G$, and so $\G$ is vertex-transitive, a contradiction. This completes the proof of Corollary~\ref{cor2p3}.

\end{document}